\documentclass[final,3p,12pt,times]{elsarticle}

\usepackage{amssymb}
\usepackage{amsthm}
\usepackage{mathptmx}
\usepackage{newtxtext,newtxmath}
\usepackage{comment}
\usepackage{wasysym}
\newdefinition{rmk}{Remark}
\DeclareMathOperator{\tr}{tr}

\begin{document}

\begin{frontmatter}


\title{On improving the numerical convergence of highly nonlinear elasticity problems}

\author[swansea]{Yue Mei}
\author[chile,chile2]{Daniel E. Hurtado}
\author[swansea]{Sanjay Pant}
\author[swansea]{Ankush Aggarwal\corref{cor1}}
\address[swansea]{Zienkiewicz Centre for Computational Engineering,\\ College of Engineering, Swansea University, Swansea SA1 8EN, UK}
\address[chile]{Department of Structural and Geotechnical Engineering, School of Engineering, \\Pontificia Universidad Cat\'{o}lica de Chile, Vicu\~{n}a Mackenna Santiago, 4860, Chile}
\address[chile2]{Institute for Biological and Medical Engineering, Schools of Engineering, Medicine and Biological Sciences, Pontificia Universidad Cat\'{o}lica de Chile, Vicu\~{n}a Mackenna Santiago, 4860, Chile}
\cortext[cor1]{Corresponding Autho r\\\textit{Email address}: a.aggarwal@swansea.ac.uk (Ankush Aggarwal)}

\begin{abstract}
Finite elasticity problems commonly include material and geometric nonlinearities and are solved using various numerical methods. However, for highly nonlinear problems, achieving convergence is relatively difficult and requires small load step sizes. In this work, we present a new method to transform the discretized governing equations so that the transformed problem has significantly reduced nonlinearity and, therefore, Newton solvers exhibit improved convergence properties. We study exponential-type nonlinearity in soft tissues and geometric nonlinearity in compression, and propose novel formulations for the two problems. We test the new formulations in several numerical examples and show significant reduction in iterations required for convergence, especially at large load steps. Notably, the proposed formulation is capable of yielding convergent solution even when 10 to 100 times larger load steps are applied. The proposed framework is generic and can be applied to other types of nonlinearities as well. 
\end{abstract}


\end{frontmatter}


\section{Introduction}
\label{intro}

With the advance of computational techniques, nonlinearities are becoming increasingly commonplace in mechanical models of solids. In finite elasticity problems, these nonlinearities can arise from different sources: material, geometry, and boundary conditions. Analytical solutions are rarely obtainable for nonlinear problems, making numerical solutions a necessity. Gradient-based methods are commonly used to numerically solve nonlinear problems, where, irrespective of the nature or degree of nonlinearity, the governing equations are linearized to obtain a Newton- or quasi-Newton-based iterative algorithm for finding the solution. For nonlinear problems, the solution is not obtained in a single step. Instead, the problem is commonly divided into smaller load steps, which are solved sequentially. High nonlinearity can lead to slow convergence, or even non-convergence, and limit the permissible load step size resulting in extremely slow computations.

A large number of studies have focused on developing preconditioners that decrease the condition number of a linear system of equations and improve their solver convergence \cite{greenbaum1997iterative}. However, those techniques usually apply to linearized systems and do not take into account specific nature of the problem. In this paper, we study the convergence properties of highly nonlinear static elasticity problems and propose a novel formulation to improve them by applying a transformation \emph{before} linearization. The proposed formulation is problem-specific, and we focus on two distinct types of nonlinearities: 1) Material nonlinearity, a common feature of constitutive models for soft biological tissues, which typically contain an exponential function \citep[chap. 4]{Maurel1998}; and 2) Geometric nonlinearity, which plays an important role in large compression problems as  the compressive force required grows rapidly with increasing compression for all material types. The new formulation is used to solve the finite-element (FE) discretization of nonlinear problems \citep{hughes2012finite, belytschko2013nonlinear}. Nonetheless, the proposed framework is equally applicable to other Galerkin-based discretization schemes, such as isogeometric analysis \citep{hughes2005isogeometric} and meshfree methods \cite{CHEN1996195}.

The mathematical foundation of the proposed method is developed in Section~\ref{analysis}. We test the performance of the new formulation using uniaxial stretching examples in Section~\ref{numerical-examples} and present the results in Section~\ref{results}. We then use the formulation to solve three practical problems and establish its improved convergence in Section~\ref{practical-problems}. Finally, we discuss the significance of the proposed framework and its equivalence to the idea of preconditioning in Section~\ref{Discussion}, before ending with a conclusion in Section~\ref{Conclusion}. 

\section{Analysis}
\label{analysis}
\subsection{Newton's method}
For finding the value of $x=x^*$ such that $g(x^*)=0$ for a general nonlinear function $g(x)$, we expand the function in Taylor's series about a point $x_n$ (guess of the solution) up to second order: 
\begin{equation}
g(x^*)=g(x_n) + \left.\frac{\partial g}{\partial x} \right|_{x_n} (x^*-x_n) + \frac{1}{2}\left.\frac{\partial^2 g}{\partial x^2} \right|_{\zeta} (x^*-x_n)^2=0,  \label{eq:Taylor}
\end{equation}
for some $\zeta \in [x_n,x^*]$. Here, remainder theorem is used to write the last term
in Eq.~\eqref{eq:Taylor}, which can also be written in other forms \citep{apostol1967calculus}. Rearranging Eq.~\eqref{eq:Taylor} we get
\begin{equation}
x^* = x_n - \frac{g(x_n)}{g'(x_n)} - \frac{g''(\zeta)}{2g'(x_n)} (x^*-x_n)^2, \label{eq:Taylor2}
\end{equation}
where $'$ denotes the differentiation with respect to the function's argument. However $\zeta$ is unknown. Therefore, neglecting the second order term gives us the classical Newton's method for determining the next solution in the iterative procedure
\begin{equation}
x_{n+1} = x_n - \frac{g(x_n)}{g'(x_n)}. \label{eq:newton}
\end{equation}
Using the definition of error at any step $e_n=|x^*-x_n|$, Eqs.~\eqref{eq:Taylor2} and \eqref{eq:newton} can be rearranged to get the evolution of error
\begin{equation}
e_{n+1} = \left|\frac{g''(\zeta)}{2g'(x_n)} \right| e_n^2. \label{eq:newton-error}
\end{equation}
Eq.~\eqref{eq:newton-error} proves the quadratic convergence of Newton's method. However, the convergence is dependent on $\left|\frac{g''(\zeta)}{2g'(x_n)} \right|$, which vanishes only for a linear function and is non-zero for any nonlinear function $g$. We define 
\begin{equation}
C(x_n,\zeta):\stackrel{\text{def}}{=}\left|\frac{g''(\zeta)}{2g'(x_n)} \right|.\label{measure-nonlinearity}
\end{equation}
We note that we know the current guess $x_n$, but value of $\zeta$ lies anywhere between $x_n$ and the solution $x^*$. Therefore, we define the maximum nonlinearity measure at the current point as 
\begin{equation}
N(x_n,x^*):\stackrel{\text{def}}{=}\sup_{\zeta\in(x_n,x^*)}\left|\frac{g''(\zeta)}{2g'(x_n)}\right|. \label{measure-nonlinearity2}
\end{equation}
For the cases where such a supremum cannot be determined, we use a local measure of nonlinearity 
\begin{equation}
\bar{C}(x_n):\stackrel{\text{def}}=C(x_n,\zeta=x_n) = \left|\frac{g''(x_n)}{2g'(x_n)}\right|. \label{measure-nonlinearity3}
\end{equation}
These defined measures will be used to quantify the degree of nonlinearity for elasticity problems defined next.

\subsection{Nonlinear elasticity problem}
\label{nonlinear-problem-section}
Given a domain $\Omega\subset\mathbb{R}^n$, a nonlinear elasticity problem involves finding a deformation mapping, i.e. a map from the reference to the deformed positions $\phi:\boldsymbol{X}\rightarrow\boldsymbol{x}$ over the domain $\Omega$, such that it satisfies the mechanical governing equations under given loading and boundary conditions. Following the standard definitions, the deformation gradient is $\mathbf{F}=\nabla_{\boldsymbol{X}}\phi=\partial \boldsymbol{x}/\partial \boldsymbol{X}$ and right Cauchy-Green deformation tensor is $\mathbf{C}=\mathbf{F}^\top\cdot\mathbf{F}$ with first three isotropic invariants 
\begin{align}
{I}_1 &= \tr \left(\mathbf{C}  \right) , \nonumber \\ 
{I}_2 &= \frac{1}{2} \left[ \tr^2 \left( \mathbf{C} \right) - \tr \left(\mathbf{C}^2  \right) \right] \text{ and} \nonumber \\
J &= \sqrt{\det \left(\mathbf{C} \right)}.
\end{align}
Green-Lagrange strain tensor is $\mathbf{E}=(\mathbf{C}-\mathbf{I})/2$ ($\mathbf{I}$ being the identity tensor). Stretch along any direction $\boldsymbol{N}$ is given by $\lambda=\sqrt{\boldsymbol{N}\cdot\mathbf{C}\boldsymbol{N}}$. The strain energy density is $W(\mathbf{F})$, from which stresses are derived via differentiation. The first Piola-Kirchhoff (PK) stress $\mathbf{P}=\partial W/\partial \mathbf{F}$, second PK stress $\mathbf{S}=\mathbf{F}^{-1}\cdot\mathbf{P}$ and Cauchy's stress $\boldsymbol{\sigma}=J^{-1}\mathbf{P}\cdot\mathbf{F}^\top$.

Adopting a discretization, the deformation mapping is approximated in terms of $N$ nodal positions using shape functions $\Psi_i$, such that $\phi^h(\boldsymbol{X}) = \sum_i\Psi_i(\boldsymbol{X})\boldsymbol{x}_i$. Discretized equations for a static  nonlinear elasticity problem can be written as
\begin{equation}
R_i(\mathbf{x})=f_i^{\text{int}}(\mathbf{x}) - f_i^{\text{ext}}(\mathbf{x}) = 0\;\forall i=1,\dots,N,\label{original-eq}
\end{equation}
to be solved for the vector with position of all nodes $\mathbf{x}$ \footnote{We differentiate $\mathbf{x}$ (the vector of all nodes' deformed positions) from $\boldsymbol{x}$ (the deformed position field)}. Here $i$ is the node number where Eq.~\eqref{original-eq} needs to be satisfied by calculating the updated displacement of all nodes $\mathbf{x}$. $R$ represents the residual; superscript $^{\text{int}}$ refers to the internal forces due to stresses and $^{\text{ext}}$ refers to the sum of all other forces -- traction forces, body forces, constraint forces, contact forces etc. In standard formulations, Eq.~\eqref{original-eq} is linearized about the current guess of the node positions $\mathbf{x}_n$, taking the form
\begin{equation}
\sum\limits_j(K^{\text{int}}_{ij}-K^{\text{ext}}_{ij})\Delta x_j = f_{i}^{\text{ext}}(\mathbf{x}_n)  - f_{i}^{\text{int}}(\mathbf{x}_n) \;\forall i=1,\dots,N.
\label{linear-eq}
\end{equation}
Here $K_{ij}$ is an element of the stiffness matrix (derivative of the forces w.r.t. node positing $x_j$) for both internal and external forces. System of equations (\ref{linear-eq}) for all the nodes is then iteratively solved. The expression for the internal forces derives from the constitutive (stress-strain) relationship of the elastic material. For total Lagrangian formulation, the nodal component of the internal force reads
\begin{equation}
f_i^{\text{int}} = 
\int_\Omega  \mathbf{P}\cdot \nabla_{\boldsymbol{X}} \Psi_i \; d\Omega,
\end{equation}
which is linear in the first PK stress, and the stretch along any axis is linear in displacement in that direction $\lambda \sim \boldsymbol{x}$. Thus, the primary nonlinearity in the internal force comes from the stress-stretch relationship.

\subsection{Proposed generalized framework}
\label{general-method}
Instead of solving the standard equation (\ref{original-eq}), we propose to solve a transformed equation
\begin{equation}
\mathcal{T} \left(f_i^{\text{int}}(\mathbf{x}) \right) = \mathcal{T} \left(f_i^{\text{ext}}(\mathbf{x}) \right)\;\forall i=1,\dots,N\label{transformed-eq}
\end{equation}
for a pre-determined bijective transformation $\mathcal{T}:\mathbb{R}\rightarrow\mathbb{R}$. 
\begin{rmk}
We note two important points about the proposed transformation (\ref{transformed-eq}):
\begin{enumerate}
\item The bijection property ensures that the solutions to Eqs.~\eqref{transformed-eq} and \eqref{original-eq} are identical.
\item Different transformations can be applied to different nodes and/or along different axes, i.e., a mixed method can be used.
\end{enumerate}
\end{rmk}

The transformation must be such that it decreases the degree of nonlinearity, which mainly comes from the stress-stretch relationship. If we use a generic scalar stress measure $\sigma$ to denote the Cauchy, 1st PK, or 2nd PK stress, then as an approximation we seek to reduce the nonlinearity of $\sigma(\lambda)$. Thus, ideally we would like to find transformation such that $\mathcal{T}(\sigma(\lambda))$ is linear. In other words, its second derivative must be zero (or as close to zero as possible):
\begin{equation}
\frac{d^2 \left[ \mathcal{T}\left( \sigma(\lambda) \right) \right]}{d\lambda^2}=\mathcal{T}''\left( \frac{d\sigma(\lambda)}{d\lambda} \right)^2 + \mathcal{T}' \frac{d^2\sigma(\lambda)}{d\lambda^2}=0.\label{eq-tau}
\end{equation}

\begin{rmk}
The exact solution of differential Eq.~(\ref{eq-tau}) is $\mathcal{T}(\sigma)= c_1\lambda(\sigma)+c_2$, \emph{i.e.} the inverse of the stress-stretch relationship. However, it may not be practical to use it because of the difficulty in obtaining its explicit mathematical expression or the computational expense of its calculation. Instead, using the main nonlinear terms in a specific stress-stretch relation, our aim is to determine a \emph{simple} transformation that \emph{reduces} the nonlinearity. We will explore this in the next section and, for now, assume that such a $\mathcal{T}$ has been determined.
\label{rmk-inverse}
\end{rmk}

Once the transformation has been determined, it is applied to the discretized force balance to obtain Eq.~(\ref{transformed-eq}). Thereupon, similar to the standard formulation, both sides are linearized about point $\mathbf{x}_n$,
\begin{equation}
\mathcal{T} \left(f_i^{\text{int}}(\mathbf{x}_n) \right) 
+ \sum\limits_j \mathcal{T}'\left(f_i^{\text{int}}(\mathbf{x}_n) \right) 
\left.\frac{\partial f_i^{\text{int}}}{\partial {x}_j} \right|_{\mathbf{x}_n} \Delta x_j = 
\mathcal{T} \left(f_i^{\text{ext}}(\mathbf{x_n}) \right) 
+ \sum\limits_j \mathcal{T}'\left(f_i^{\text{ext}}(\mathbf{x}_n) \right)
\left.\frac{\partial f_i^{\text{ext}}}{\partial {x}_j} \right|_{\mathbf{x}_n} \Delta x_j.
\label{transformed-eq2}
\end{equation}
Using the definition of the stiffness matrix and rearranging we get
\begin{equation}
\sum\limits_j\left[ \mathcal{T}'(f_i^{\text{int}}(\mathbf{x}_n))K^{\text{int}}_{ij} - \mathcal{T}'(f_i^{\text{ext}}(\mathbf{x}_n))K^{\text{ext}}_{ij} \right] \Delta x_j= \mathcal{T} \left(f_i^{\text{ext}}(\mathbf{x_n}) \right) -\mathcal{T} \left(f_i^{\text{int}}(\mathbf{x}_n) \right) ,\label{transformed-eq3}
\end{equation}
Defining a new symbol $\mu_i=\mathcal{T}'(f_i^{\text{ext}})/\mathcal{T}'(f_i^{\text{int}})$, we can write Eq.~\eqref{transformed-eq3} as
\begin{equation}
\sum\limits_j\mathcal{T}'(f_i^{\text{int}}(\mathbf{x}_n))\left(K^{\text{int}}_{ij}  - \mu_i (\mathbf{x}_n) K^{\text{ext}}_{ij} \right) \Delta x_j =  
\mathcal{T} \left(f_i^{\text{ext}}(\mathbf{x}_n) \right) - \mathcal{T} \left(f_i^{\text{int}}(\mathbf{x}_n) \right), \label{transformed-linear-eq}
\end{equation}
Furthermore, using the approximation that at $\mathbf{x}_n$ the internal and external forces are of similar magnitude (i.e. $f_i^{\text{int}}(\mathbf{x}_n) \approx f_i^{\text{ext}}(\mathbf{x}_n) \Rightarrow \mu_i(\mathbf{x}_n)\approx1$), we obtain
\begin{equation}
\sum\limits_j\mathcal{T}'(f_i^{\text{int}}(\mathbf{x}_n))\left(K^{\text{int}}_{ij}  -  K^{\text{ext}}_{ij} \right) \Delta x_j =  
\mathcal{T} \left(f_i^{\text{ext}}(\mathbf{x}_n) \right) - \mathcal{T} \left(f_i^{\text{int}}(\mathbf{x}_n) \right). \label{transformed-linear-eq2}
\end{equation}
 Eq.~\eqref{transformed-linear-eq2} can also be written as
\begin{equation}
\sum\limits_j\left(K^{\text{int}}_{ij}  -  K^{\text{ext}}_{ij} \right) \Delta x_j =  
\frac{\mathcal{T} \left(f_i^{\text{ext}}(\mathbf{x}_n) \right) - \mathcal{T} \left(f_i^{\text{int}}(\mathbf{x}_n) \right)}{\mathcal{T}' \left(f_i^{\text{int}}(\mathbf{x}_n) \right)}. \label{transformed-linear-eq3}
\end{equation}

\begin{rmk}
In Eq.~(\ref{transformed-linear-eq3}), if we expand the $\mathcal{T} \left(f_i^{\text{int}}(\mathbf{x}_n) \right)$ term on the right hand side using Taylor's series about $f_i^{\text{ext}}(\mathbf{x}_n)$ and then truncate to the first order, assuming that the difference $f_i^{\text{int}}(\mathbf{x}_n)-f_i^{\text{ext}}(\mathbf{x}_n)$ is small, we get the standard formulation (\ref{linear-eq}) back. Thus, the transformed equation is equivalent to the standard equation if the load step is small enough, but it becomes increasingly different as we increase the load step size.
\end{rmk}

\begin{rmk}
Although Eq.~(\ref{transformed-eq3}) could be implemented as is, the approximation $\mu_i(\mathbf{x}_n)\approx1$ provides a highly simplified form. Comparing with the standard formulation (\ref{linear-eq}), we note that (\ref{transformed-linear-eq3}) does not involve changing the stiffness matrix, and only the right hand side is different. Therefore, the modification required for the proposed formulation is minimal at the solver stage. 
\end{rmk}

\begin{rmk}
The difference between Eqs. (\ref{transformed-linear-eq2}) and (\ref{transformed-linear-eq3}) is akin to linear preconditioning. However, it is not clear which of the two will have a lower condition number. Therefore, we assume that the condition number of the standard stiffness matrix is low enough and use (\ref{transformed-linear-eq3}).
\end{rmk}

To determine the form of transformation $\mathcal{T}$ and whether/when it is useful to use the transformed equation~(\ref{transformed-eq}) instead of the standard equation~(\ref{original-eq}), we next look at specific nonlinearities.


\subsection{Material non-linearity}
\label{material-non-linearity-section}
Many of the constitutive models for soft tissues contain an exponential function, which is the primary source of nonlinearity. Therefore, we assume a highly-simplified form $\sigma\sim\exp(\lambda)$, for which one could reduce the nonlinearity by taking a logarithm, i.e., $\mathcal{T}\equiv\log$. We use Eq.~\eqref{transformed-linear-eq3} for the nodes where the log transformation is applicable, and keep Eq.~\eqref{linear-eq} for other nodes. In general, we write the linearized system of equations as
\begin{equation}
\sum\limits_j K_{ij} \Delta x_j = \bar{R}_i (\mathbf{x}_n),
\label{modified-equation}
\end{equation}
where the modified residual
\begin{equation}
\bar{R}_i (\mathbf{x}_n) = \left\{ \begin{array}{cc}
f_{i}{^\text{int}}(\mathbf{x}_n) \log\left(\frac{f_{i}^{\text{ext}}(\mathbf{x}_n)}{f_{i}^{\text{int}} (\mathbf{x}_n)} \right) & \text{if Condition (\ref{eq:log-condition}) is satisfied} \\
f_{i}^{\text{ext}}(\mathbf{x}_n)  - f_{i}^{\text{int}}(\mathbf{x}_n) & \text{otherwise}
\end{array}\right.
\label{residual-extension}
\end{equation}
and the Condition, for some tolerance $\mathrm{TOL}$, is
\begin{equation}
|f_{i}^{\text{ext}}(\mathbf{x}_n)| > \mathrm{TOL} \text { and }  |f_{i}^{\text{int}}(\mathbf{x}_n)| > \mathrm{TOL} \text{ and } \frac{f_{i}^{\text{ext}}(\mathbf{x}_n)}{f_{i}^{\text{int}}(\mathbf{x}_n)} > 0.
\label{eq:log-condition}
\end{equation}
This is satisfied when both internal and external forces at the current positions $\mathbf{x}_n$ are non-zero and of the same sign. We call this ``log formulation'' that deals with the exponential nonlinearity.

\subsubsection{Error comparison}

To determine whether the transformed formulation leads to an advantage over the standard formulation, we look at the measure of nonlinearity $N$, Eq.~(\ref{measure-nonlinearity2}), for functions with a single unknown. We start with a simple exponential function: $g(x): Ae^{Bx}=H$ to be solved for a given constant $H$, such as that used to derive the log formulation. The standard formulation gives 
\begin{equation}
N_\text{standard} = \sup_{\zeta\in(x_n,x)}\left|Be^{B(\zeta-x_n)} \right|=Be^{B(x-x_n)}=Be^{B\Delta x}.
\end{equation}
It is evident that the error in the Newton's method will increase as we increase $\Delta x$, \emph{i.e.} if the initial guess is farther away from the solution. In some cases, Newton's method may not even converge. Whereas, applying the transformation $\log(Ae^{Bx})=\log(H)$, we get a linear equation and $N_\text{transform}=0$ identically. Thus, the transformed method will always converge in a single iteration, and it will always be better to use the new formulation compared to the standard formulation. 

A more realistic model of force-displacement relation is $g(x):A(e^{Bx}-1)=H$ as the left hand side is zero at $x=0$. In this case, $N_\text{standard}=Be^{B\Delta x}$ for the standard formulation. On the other hand, using the transformation $\log\left[A(e^{Bx}-1)\right]=\log[H]$, we get
\begin{equation}
N_\text{transform} = \sup_{\zeta\in(x_n,x)}\left|\frac{-\left(\frac{1}{f(\zeta)}ABe^{B\zeta}\right)^2 + \frac{1}{f(\zeta)}AB^2e^{B\zeta}}{\frac{1}{f(x_n)}ABe^{Bx_n}} \right|=\frac{B}{e^{Bx_n}-1}.
\end{equation}
Thus, the nonlinearity in the new formulation is non-zero but does not depend on $\Delta x$. It only depends on the starting guess $x_n$ and decreases as we increase $x_n$. Therefore, a \emph{conservative condition} for the new formulation to perform better than the standard formulation is $N_\text{standard}> N_\text{transform}$, i.e.,
\begin{equation}
\frac{Be^{Bx}}{e^{Bx_n}} > \frac{B}{e^{Bx_n}-1}.
\end{equation}
After rearranging, we can write the above condition as
\begin{equation}
H=A(e^{Bx}-1) > \frac{A}{e^{Bx_n}-1}.
\end{equation}
Thus, unless $Bx_n$ is extremely small, for $H\gtrapprox A$, the new formulation is expected to outperform the standard formulation. 


\subsection{Geometric Non-linearity}
\label{geometric-non-linearity}
In compression, elastic solids exhibit geometric nonlinearity, even for linear constitutive models. Thus, we look at the uniaxial compression case with stretch $\lambda<1$ in the compression direction and a volume-preserving material. The associated deformation gradient is $\mathbf{F}=\text{diag}\left[ \lambda, 1/\sqrt{\lambda},1/\sqrt{\lambda} \right]$ for compression along the first axis. Under this deformation, stress along the first axis for any incompressible isotropic material is given by
\begin{equation}
\sigma^{(n)}(\lambda)= 2\frac{\partial W}{\partial I_1}\frac{1}{\lambda^n}\left(\lambda^2 -\frac{1}{\lambda} \right) + 2\frac{\partial W}{\partial I_2}\frac{1}{\lambda^n}\left(\lambda -\frac{1}{\lambda^2} \right),
\end{equation}
where $n=0$, 1 and 2 correspond to Cauchy, 1st PK and 2nd PK stresses, respectively. 

As a simplification, we assume $\partial W/\partial I_1=1$ and $\partial W/\partial I_2=0$ (i.e. a neo-Hookean constitutive model with unit shear modulus), and thus 
\begin{equation}
\sigma^{(n)}(\lambda)=\lambda^{-n}(\lambda^2-1/\lambda).\label{eq:all-stresses}
\end{equation}
It is difficult to find an inverse function of the above relation (see Remark \ref{rmk-inverse}). Substituting the above stress-stretch relation into Eq.~(\ref{eq-tau}) and rearranging, we get
\begin{equation}
{\mathcal{T}^{\prime \prime }}\left( {{{\left( {{\sigma ^{(n)}}} \right)}^2} + {\beta ^2}\left( \lambda  \right)} \right) + {\mathcal{T}^\prime }\left( {2{\sigma ^{(n)}} + \gamma \left( \lambda  \right)} \right) = 0, \label{eq-tau-comp}
\end{equation}
where 
\begin{align}
\beta ^2(\lambda) &= \frac{n}{n+2}\lambda^{-2n-2} - 2\frac{n^2-5n-6}{(n+1)(n+2)}\lambda^{1-2n} + \frac{n^2-11n+6}{(n+1)(n+2)}\lambda^{4-2n} \text{ and} \nonumber \\
\gamma \left( \lambda  \right) &= \frac{ - 12n }{(n + 1)(n + 2)} \lambda^{2 - n}.
\end{align}
We note that $\beta(\lambda)$ is unbounded for 1st and 2nd PK stresses as $\lambda\rightarrow 0$, making it difficult to find a closed form solution of Eq.~\eqref{eq-tau-comp} for $n=$1 or 2. However, for $n=0$, i.e. Cauchy stress, we have
\begin{equation}
{\cal T}''\left[ \left(\sigma^{(0)} \right)^2 + {\beta ^2}\left( \lambda  \right) \right] + 2{\cal T}'\sigma^{(0)} = 0,\label{eq-Cauchy-T}
\end{equation}
where $\beta^2(\lambda)=3\lambda^4+6\lambda$ varies between 9 and 0 for compression. If we approximate $\beta^2(\lambda)$ as a constant, the solution to  Eq.~(\ref{eq-Cauchy-T}) is
\begin{equation}
\mathcal{T}(\sigma) = \frac{c_1}{\beta} \tan^{-1}\left(\frac{\sigma}{\beta}\right) + c_2.
\end{equation}
We note that the values of the integration constants $c_1$ and $c_2$ do not affect the transformation as long as $c_1\ne0$. Therefore, we arbitrarily choose $c_1=\beta$ and $c_2=0$ and arrive at the transformation $\mathcal{T}(\sigma) = \tan^{-1}\left(\sigma/\beta\right)$. 

\begin{figure}[!t]
  \centering
  \includegraphics[width=0.9\textwidth]{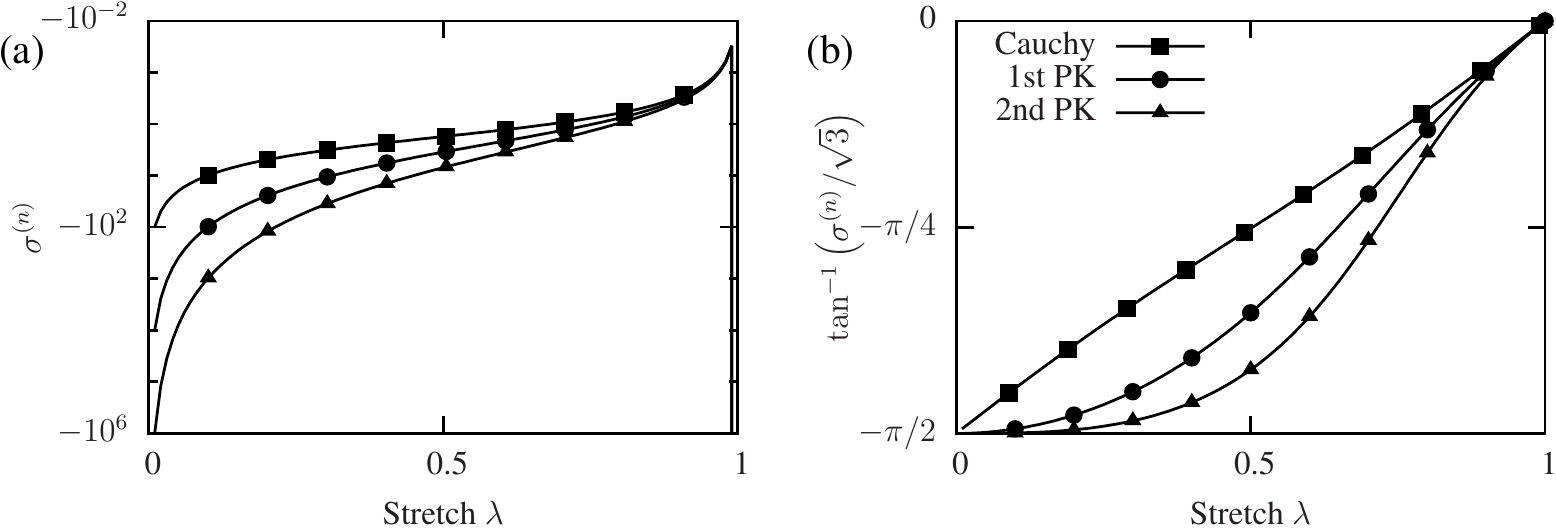}
  \caption{(a) Stress variation with respect to stretch $\lambda$ for a neo-Hookean constitutive
model with unit shear modulus under uniaxial compression; (b) the corresponding transformed stresses remain bounded and show a significantly reduced nonlinearity\label{arctan-stress-fig}}
\end{figure}

However, the value of $\beta$ remains undetermined. We use an approximate average value of $\beta^2(\lambda)$ in $\lambda\in(0,1)$: $\langle\beta^2\rangle  = \int\limits_0^1 {\beta^2 ( \lambda  )\, d\lambda } = \int\limits_0^1 3\lambda^4+6\lambda\, d\lambda  \approx 3$. The Cauchy stress becomes nearly-linear after applying this transformation (Fig.~\ref{arctan-stress-fig}), as well as the transformed 1st and 2nd PK stresses show a significantly reduced nonlinearity.

The internal force $f^\text{int}\propto\sigma$ up to an unknown multiplicative constant $\gamma$, which is a function of the shape functions and mesh density: $f^\text{int}\sim\gamma\sigma$. Therefore, the transformation that linearizes the stress, may not linearize the internal force. For material nonlinearity of the exponential type, this is not an issue since we get $\log(\gamma \sigma)= \log(\gamma) + \log(\sigma)$, and the multiplicative constant factors out resulting in a linear force-displacement relation. However, in general, even if $\mathcal{T}(\sigma)$ linearizes the stress-stretch relation, it does not imply that $\mathcal{T}(f^\text{int})$ will also be linear, such as for $\mathcal{T}\equiv\tan^{-1}$. In order to resolve this, we introduce a finite, non-zero factor $\alpha$ into our transformation
\begin{equation}
\mathcal{T}(f^\text{int})=\tan^{-1}\left(\alpha f^\text{int} \right).
\end{equation}
To compute $\alpha$, we use the fact that $f^\text{int}(\lambda=1)=0$ and $f^\text{int}(\lambda=0)=\infty$. As a result, $\tan^{-1}\left(\alpha f^\text{int} \right)=0$ and $\pi/2$ at $\lambda=1$ and 0, respectively. Furthermore, for $\lambda\in(0,1)$, we want $\tan^{-1}\left(\alpha f^\text{int} \right)$ to be linear in $\lambda$. Therefore, if at a stretch value of $\lambda_i$, the internal force $f_i^\text{int}$ is known (for example at current iteration), we calculate $\alpha$ at node $i$ using the relation
\begin{equation}
\tan^{-1} \left( \alpha_i f_i^\text{int}(\lambda_i) \right) =  \frac{\pi}{2} \left( 1- \lambda _i \right),
\end{equation}
and update its value at every iteration. In general, we again use Eq.~\eqref{modified-equation} with the modified residual
\begin{equation}
\bar{R}_i (\mathbf{x}_n) = \left\{ \begin{array}{cc}
\frac{1+ \left(\alpha_if_{i}{^\text{int}}\right)^2}{\alpha_i}  \left(\tan^{-1} (\alpha_if_{i}^{\text{ext}})- \tan^{-1}(\alpha_if_{i}^{\text{int}}) \right) & \text{if Condition (\ref{eq:tan-condition}) is satisfied} \\
f_{i}^{\text{ext}}  - f_{i}^{\text{int}} & \text{otherwise}
\end{array}\right..
\label{residual-compression}
\end{equation}
Here, the forces are calculated at $\mathbf{x}_n$ and the Condition, for some tolerance $\mathrm{TOL}$, is
\begin{equation}
|f_{i}^{\text{ext}}(\mathbf{x}_n)| > \mathrm{TOL} \text { and } \alpha_i \text{ is computed from previous iteration}.
\label{eq:tan-condition}
\end{equation}
We call this ``arctan formulation'' that reduces the geometric nonlinearity.

\begin{rmk}
In FEM, the forces are calculated at the nodes whereas strains and stretches are calculated at the integration or Gauss points of the elements. However, we only need an approximate value of $\lambda$ at the nodes. Therefore, one may average the stretches from connected elements, calculate stretch at the node, and use it in the above relation.
\end{rmk}

\subsubsection{Error Comparison}
\begin{figure}[!t]
 \centering
 \includegraphics[width=0.9\textwidth]{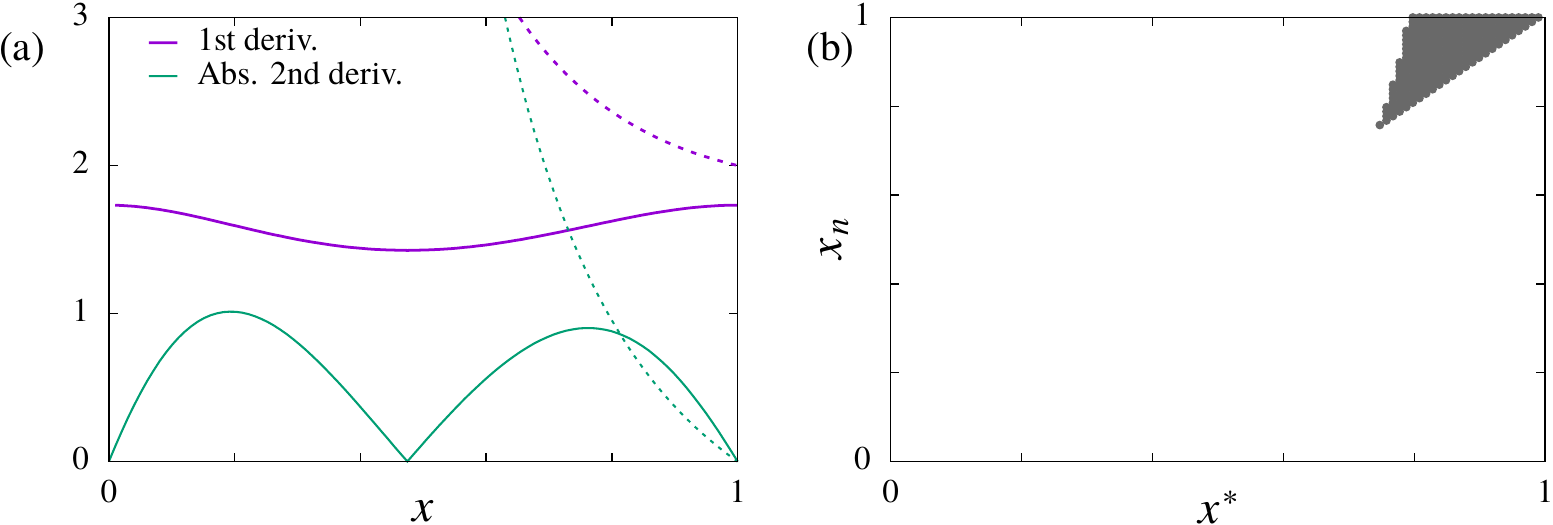}
 \caption{(a) Derivatives of the stress in compression before (dashed line) and after the arctan transformation (solid line); (b) Plot of the area where $N_{\text{standard}} < N_{\text{transform}}$ denoting the region where standard formulation \emph{may} be better than the new formulation\label{atan-error}}
\end{figure}

In order to compare the nonlinearity before and after transformation, we look at the Cauchy stress used to derive the arctan formulation (Eq.~\ref{eq:all-stresses} with $n=0$). Accordingly, we construct a function $g(x):x^2-1/x=H$ for a given constant $H$. Its first and second derivatives go to infinity as $x\rightarrow 0$, whereas after transformation $\tan^{-1}\left(\frac{x^2-1/x}{\sqrt{3}}\right)=\tan^{-1}\left(\frac{H}{\sqrt{3}}\right)$ the derivatives remain bounded in $x\in(0,1)$ (Fig.~\ref{atan-error}). Because of the oscillatory nature of derivatives, an analytical comparison of nonlinearity measure $N$ is not possible for the two cases. Instead, we perform a numerical comparison and plot the region where $N_\text{standard}< N_\text{transform}$ (Fig.~\ref{atan-error}). It is clear that only in small compression cases ($x\not\ll1$), the standard formulation \emph{may} perform better than the new formulation. 

\section{Numerical Examples}
\label{numerical-examples}
To test the feasibility of the proposed formulation and its effect on the convergence, we first solve simple extension/compression problems along an axis with varying material properties and loading steps.

\subsection{Problem description}

\begin{figure}[!t]
  \centering
  \includegraphics[width=1.0\textwidth]{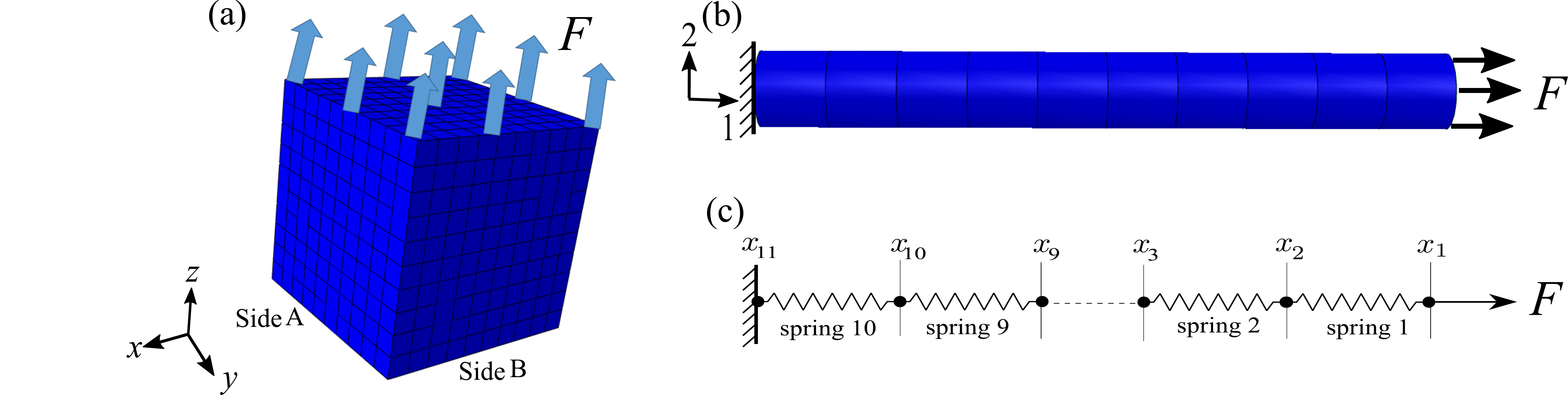}
  \caption{Schematic of the uniaxial extension/compression for (a) a three-dimensional solid , (b) axisymmetric case , and (c) a one-dimensional case  \label{schematic}}
\end{figure} 

We investigate the problem of pure uniaxial extension/compression of an isotropic hyperelastic solid and its reduction to axisymmetric and one-dimensional cases (Fig.~\ref{schematic}). In the three-dimensional (3D) case, we consider a solid cube of a unit edge length and mesh it uniformly into $10\times10\times10$ trilinear hexahedral elements. A uniform pressure loading $F$ is applied on the top face (either tensile or compressive), and the motion of the bottom face is restricted in the $z$-direction. Moreover, to remove the rigid body rotation modes, we restrict the motion of sides A and B along the $x$ and $y$ direction, respectively (Fig.~\ref{schematic}a). 

Secondly, we consider the case of an axisymmetric cylinder under uniaxial loading, so that, by symmetry, the deformation gradient $\mathbf{F}=\text{diag}\left[ \lambda_1, \lambda_2,\lambda_2 \right]$ for axial stretch $\lambda_1$ and lateral stretch $\lambda_2$. This problem is solved using one-dimensional linear elements with two degrees of freedom (DOFs) per node: axial and lateral displacements (Fig.~\ref{schematic}b). We restrict the axial motion of the left end of the cylinder, apply a uniform pressure load $F$ on its right face, and set lateral stress to be zero everywhere. Even though this problem is identical to the 3D case, its computational implementation is simpler because of the fewer degrees of freedom and only one traction boundary node. Thus, a comparison between the 3D and axisymmetric cases would allow us to study the effect of dimensionality and mesh refinement on the proposed formulation.

Lastly, the above problem is further simplified under the assumption of incompressibility, such that $\lambda_1=\lambda$ and $\lambda_2=1/\sqrt{\lambda}$. This case is also discretized into ten one-dimensional linear elements with one DOF per node (Fig.~\ref{schematic}c). We assume that the reference length of each spring is $l$, left end of the springs system is fixed, and a pressure load $F$ is applied on the right end. Comparing this case with the compressible axisymmetric case will allow us to study the effect of compressibility. In all three cases, the pressure force acts on the reference configuration and is not a follower load. That is, the effective force does not change with deformation and $\mathbf{K}^{\text{ext}}=\mathbf{0}$. 

\subsubsection{Constitutive models}
\begin{figure}[!tbp]
  \centering
  \includegraphics[width=0.85\textwidth]{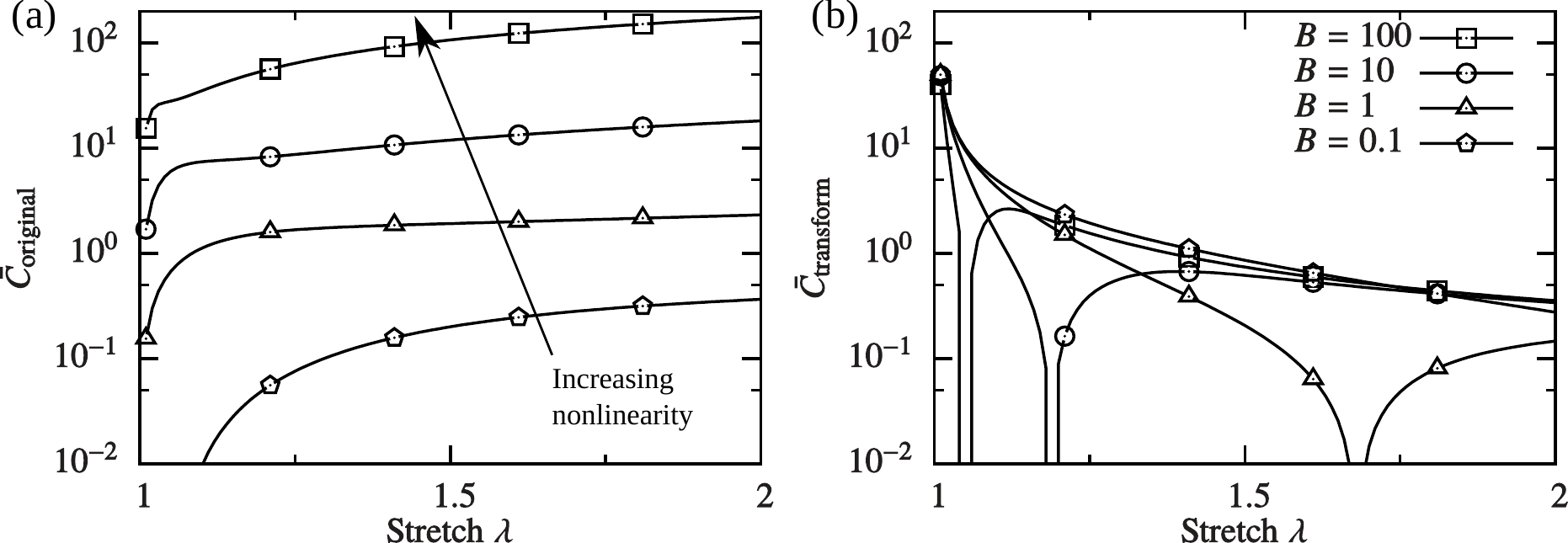}
  \caption{Local measure of nonlinearity $C$ (Eq.~\ref{measure-nonlinearity3}) of the VW model (37) for the (a) original stresses and (b) the transformed stresses for varying exponent parameter $B$\label{VMerror-f}}
\end{figure}

We consider different models for the stress-strain relationship to investigate the effect of proposed formulation in a general setting. For material nonlinearity with an exponential-type behavior, we use the Veronda-Westmann (VW) model, which is commonly used for the biomechanical response of biological tissues \citep{safshekan2016mechanical,girard2015engineering} and defines the strain energy density as
\begin{equation}
W(\mathbf{F})  = \frac{A}{B}\left[ e^{B\left( J^{-2/3}{I}_1 - 3 \right)} - 1 \right] - \frac{A}{2} \left( J^{-4/3}{I}_2 - 3 \right) + \frac{K}{2}\left( \ln J \right)^2.
\end{equation}
Here $K$ is the bulk modulus, and $A$ and $B$ are stiffness parameters. $A$ is the initial shear modulus and has the units of stress, while the nonlinearity depends on the dimensionless parameter $B$. In case of uniaxial stretch under incompressibility constraint ($J=1$), the first PK stress along stretch direction reduces to
\begin{equation}
P\left( \lambda  \right) = 2A\left( \lambda  - \frac{1}{\lambda ^2} \right) e^{B\left(\lambda ^2 + \frac{2}{\lambda } - 3\right)} - A\left( 1 - \frac{1}{\lambda ^3} \right).
\label{1stpk-eq}
\end{equation}
We note that even though this function is significantly more involved than an isolated exponential function we used to determine the transformation $\mathcal{T}\equiv\log$, for small extension, the primary nonlinearity comes from the exponential function. Hence we use log formulation to solve uniaxial extension with this material model. In order to quantify the nonlinearity of this model \eqref{1stpk-eq}, we look at the local measure of nonlinearity $\bar{C}$ (Eq.~\ref{measure-nonlinearity3}) for varying exponent parameter $B$ (Fig.~\ref{VMerror-f}). Clearly the nonlinearity of $P(\lambda)$ rapidly increases for larger values of $B$. Moreover, for a given value of $B$, the nonlinearity slightly increases with the stretch. On the other hand, the nonlinearity of the transformed stress $\log[P(\lambda)]$ decreases quickly with stretch. For the most part, the nonlinearity of transformed stress remains around unity, except sharply dipping at certain stretch values. At these points, $\bar C_{\text{transform}}$ goes to zero, i.e. the transformed stress is locally exactly linear. Comparing the nonlinearity of original and transformed stresses, we notice that, for this one-dimensional incompressible stress, the transformation reduces the nonlinearity for $B\ge1$.

\begin{figure}[!t]
  \centering
  \includegraphics[width=0.9\textwidth]{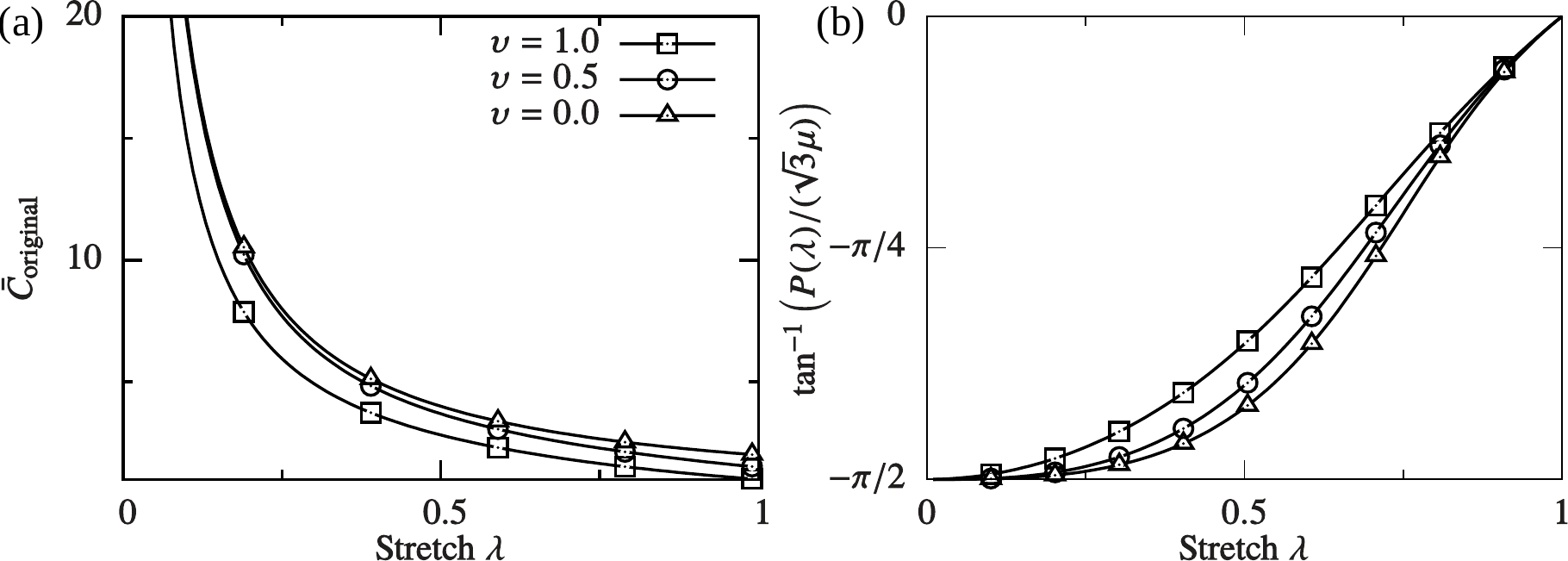}
  \caption{For the Mooney-Rivlin model \eqref{MR-1stPK} under compression and varying $\upsilon$, the measure of nonlinearity using (a) the original increases indefinitely as $\lambda\rightarrow 0$; (b) whereas the transformed 1st PK stress remains finite and close to linear \label{mooney-fig}}
\end{figure}

For geometric nonlinearity case, we study the compressible Mooney-Rivlin (MR) model, which defines the strain energy density as \citep[chap. 6]{holzapfel2000nonlinear}
\begin{equation}
W({\mathbf{F}}) = \frac{\mu}{2}\left[ \upsilon \left(J^{- 2/3}{I_1} - 3 \right) + \left( 1 - \upsilon \right)\left( J^{ - 4/3}{I_2} - 3 \right) \right] + \frac{K}{2}\left( \ln J \right)^2.
\end{equation}
Here $\mu$ is the effective shear modulus, $K$ is the bulk modulus, and $\upsilon\in[0,1]$ is a dimensionless material parameter. In this case, the uniaxial stretch and incompressibility constraint imply that the first PK stress 
\begin{equation}
P(\lambda ) = \mu \left[ \upsilon \left( \lambda  - \frac{1}{\lambda ^2} \right) + (1 - \upsilon ) \left( 1 - \frac{1}{\lambda ^3} \right) \right].\label{MR-1stPK}
\end{equation}
Substituting $\upsilon=1$ reduces the above model to the neo-Hookean model, for which the transformation $\mathcal{T}\equiv\tan^{-1}$ was determined. We vary $\upsilon$ to determine how well this transformation performs for models that deviate from neo-Hookean. The nonlinearity of the original stress increases slightly as we decrease $\upsilon$ (Fig.~\ref{mooney-fig}a). Plotting the transformed stress  $\tan^{-1}(P/(\sqrt{3}\mu))$, we notice that its nonlinearity increases as we decrease $\upsilon$ (Fig.~\ref{mooney-fig}b). Using the value of $\alpha$ based on the previous iteration will further decrease its nonlinearity.    
 
\subsubsection{Solution details}
\label{solution-section}
We use the general Eq.~\eqref{modified-equation} with the modified residual (Eq.~\eqref{residual-extension} for the VW model and Eq.~\eqref{residual-compression} for the MR model) for varying parameter values. We note that the transformation is applicable only at the traction boundary nodes, and we keep the standard formulation at the interior nodes. Thereby, we solve each problem using two load steps. We always calculate the first load step, denoted by $F_1$, using the standard formulation (i.e. without any transformation). This allows us to have a non-zero internal force for the log formulation and calculate $\alpha_i$ for the arctan formulation. Once these conditions are satisfied, we solve the transformed system of equations for the second load step under an external force $F_2$. Usually, more than two load steps are used to solve a nonlinear problem. Here, we only use two load steps as an extreme case to test convergence.

Newton's method is used to iteratively solve the discretized force balance equations. A displacement-based convergence criteria is used for all cases: $\left\| \mathbf{\Delta U} \right\|/\left\| \mathbf{U_k} \right\|<10^{-3}$, where $\left\| \mathbf{U_k} \right\|$ and $\left\| \Delta\mathbf{U} \right\|$ are $L_2$-norms of the nodal displacement vector at the current iteration and the displacement increment, respectively. Furthermore, iterations are terminated if a negative Jacobian $J$ or negative stretch ratio is detected, and we set the maximum number of iterations to 100, which is assumed to be a non-converged result. The convergence is quantified by the number of iterations taken for the solution to converge, which is computed at varying first and second load steps. For the 3D case, the proposed formulation was implemented into the open source finite element code FEBio \cite{maas2012febio}. For the other two nonlinear problems, a finite element code was implemented in Python to solve them using both standard and new formulations. 

For the VW model, we vary its nonlinearity from small ($B=1$) to large ($B=100$) and non-dimensionalize the applied traction using parameter $A$. For the compressible case, we also vary the bulk modulus $K$ and test the proposed formulation for varying degree of compressibility. For the MR model, we vary the parameter $\upsilon$ from $0$ to $1$ and non-dimensionalize the applied traction using the effective shear modulus $\mu$. We start with the simplest case of one dimensional incompressible stretch, and then solve the (compressible) axisymmetric and 3D cases (i.e.~reverse of the order presented in Fig.~\ref{schematic}). 

\section{Results}
\label{results}
\subsection{Material nonlinearity}

\subsubsection{Incompressible extension}
\label{1D-log-section}

\begin{figure}[!tbp]
  \centering
  \includegraphics[width=0.9\textwidth]{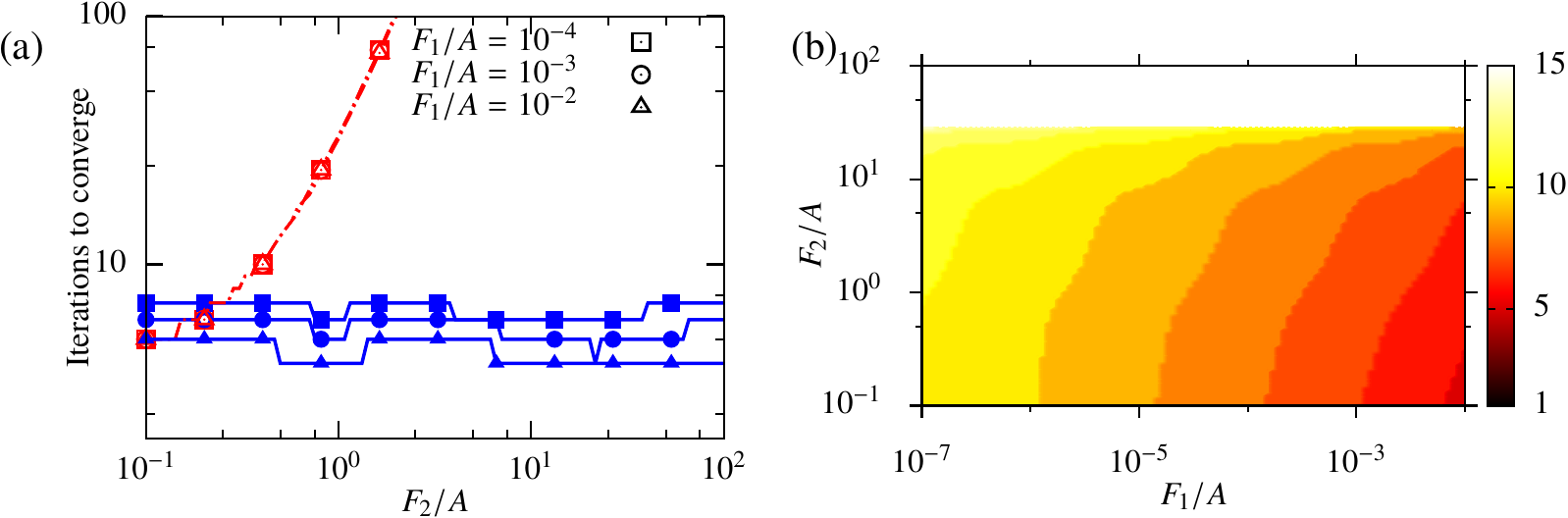}
  \caption{Effect of the starting point on iterations taken to converge with the Veronda-Westmann model: (a) for the 1D extension problem using the standard (red dashed line and open symbols) and log formulation (blue solid line and filled symbols) with $B=100$ and varying $F_1/A$, and (b) for the axisymmetric compressible extension using log formulation. \label{VM-f1}}
\end{figure}

For the axisymmetric impressible extension case, for high nonlinearity ($B=100$), the standard residual formulation takes an increasing number of iterations to converge as we increase the second load step (Fig.~\ref{VM-f1}a), and does not converge for $F_2/A\gtrapprox2$. In contrast, the log formulation converges in less than ten iterations for load step up to $F_2/A=100$. The convergence of the standard formulation is largely insensitive to the first load step, whereas for the new formulation the number of iterations decreases further as we increase the first load step (Fig.~\ref{VM-f1}b). 

\begin{figure}[!tbp]
  \centering
  \includegraphics[width=0.9\textwidth]{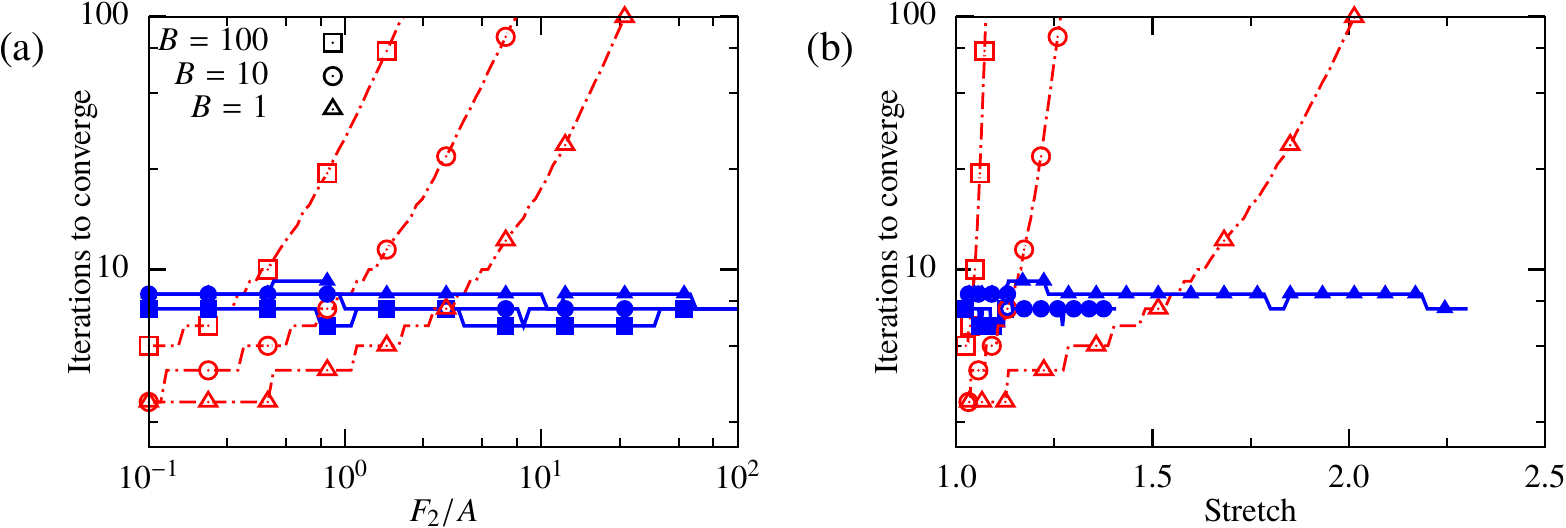}
  \caption{Iterations taken to converge for the 1D extension problem with the Veronda-Westmann model using the standard (red dashed line and open symbols) and log formulation (blue solid line and filled symbols) with $F_1/A=10^{-4}$ and varying $B$. \label{VM-f2}}
\end{figure}

If we decrease the degree of nonlinearity, the convergence of the standard formulation improves; the load step at which the method does not converge increases exponentially (Fig.~\ref{VM-f2}a). However, the number of iterations required using the log formulation remains largely constant (approximately eight), even for small nonlinearity. For small enough load step, the log formulation may take more iterations than the standard formulation. However, this disadvantage of using the log formulation is insignificant (<5 iterations difference) and exists only at small load steps. The log formulation converges at loads as high as 100 times the maximum load step possible using the standard formulation. If the number of iterations are plotted with respect to the induced stretch, a similar trend is observed (Fig.~\ref{VM-f2}b). For high nonlinearity cases, the induced stretch is limited below 1.5 as the material shows a strain-limiting behavior.

\subsubsection{Axisymmetric extension}
\begin{figure}[!tbp]
  \centering
  \includegraphics[width=0.9\textwidth]{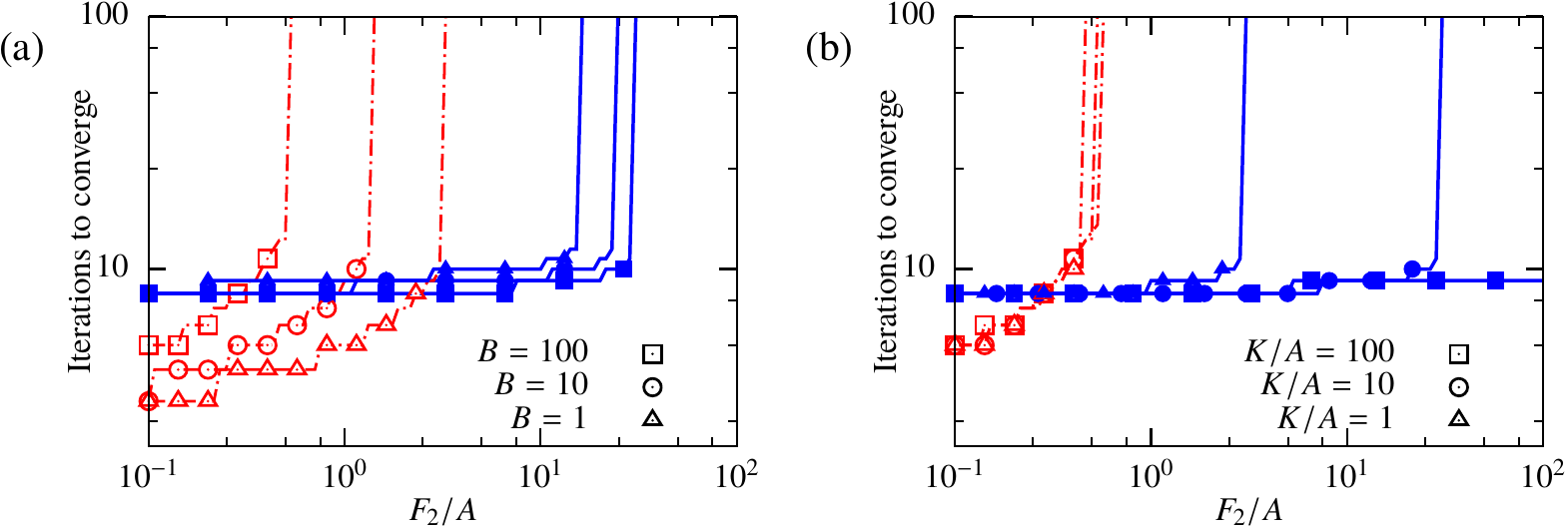}
  \caption{Iterations taken to converge for the axisymmetric extension problem with Veronda-Westmann model using the standard formulation (red dashed line and open symbols) and log formulation (blue solid line and filled symbols): (a) different values of $B$ and (b) different values of $K/A$.
   \label{VM-3D-fig}}
\end{figure}

For the compressible axisymmetric extension ($K/A=10$), we find that the standard formulation diverges for even smaller load steps (Fig.~\ref{VM-3D-fig}a). In contrast, the log formulation converges in approximately 10 iterations for all values of $B$ considered in this study (Fig.~\ref{VM-3D-fig}a), even while using 10-100 times larger load step as compared to the standard formulation. Furthermore, for large nonlinearity ($B=100$), increasing the compressibility does not significantly affect the convergence of the standard formulation (Fig.~\ref{VM-3D-fig}b). However, convergence of the proposed formulation improves as we increase the bulk modulus, as the iterations diverge at a larger load step. Similar to the incompressible case, increasing the first load step further improves the convergence of log formulation (Fig.~\ref{VM-f1}).


\subsubsection{3D uniaxial extension}
\label{3D-VM-Cube-section}
\begin{figure}[!tbp]
  \centering
  \includegraphics[width=0.9\textwidth]{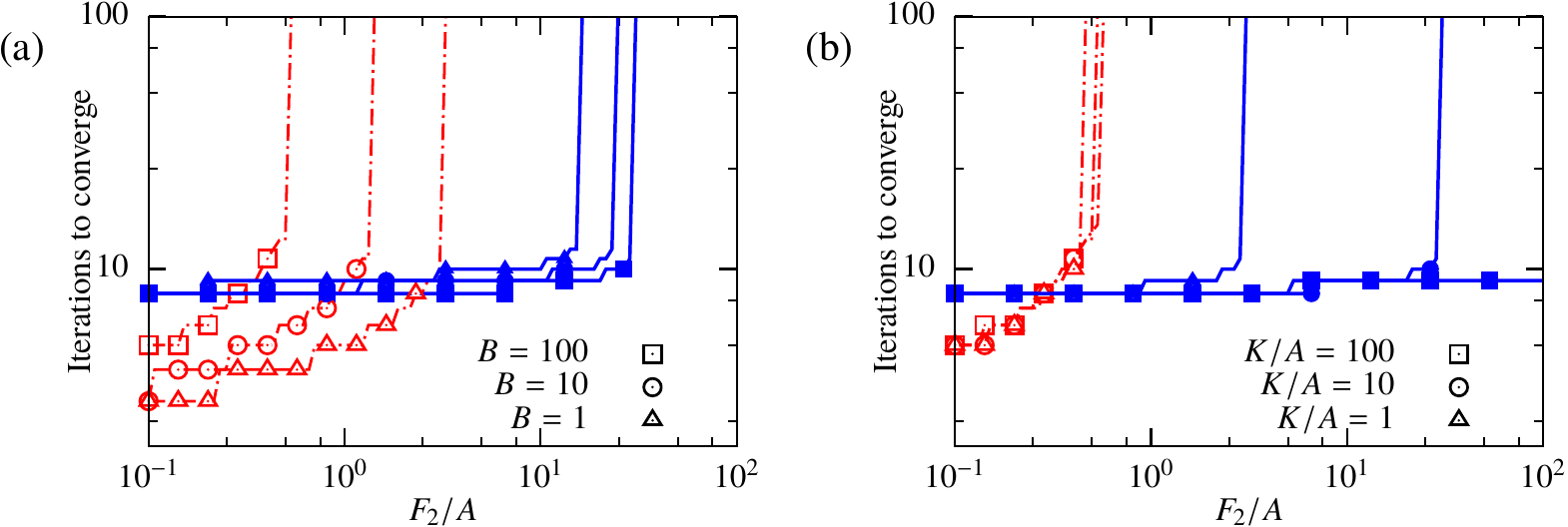}
  \caption{Iterations taken to converge for the 3D uniaxial extension problem with Veronda-Westmann model using the standard formulation (red dashed line and open symbols) and log formulation (blue solid line and filled symbols) for (a) varying $B$ and (b) varying $K/A$.\label{VM-cube-3D-fig}}
\end{figure}

Next, we relax the axisymmetric assumption and solve a three-dimensional uniaxial extension for bulk modulus $K/A=10$ (Fig.~\ref{VM-cube-3D-fig}a) and $B=100$ (Fig.~\ref{VM-cube-3D-fig}b). In spite of the higher number of DOFs compared to the previous case, the resulting number of iterations show exactly the same trend (Fig.~\ref{VM-cube-3D-fig}). Thus, the log formulation improves the convergence of a general three-dimensional elastic extension, and that improvement becomes more pronounced for nearly incompressible problem.

\subsection{Geometric nonlinearity}
\label{arctan-section}

\subsubsection{Incompressible compression}

\begin{figure}[!tbp]
  \centering
    \includegraphics[width=0.9\textwidth]{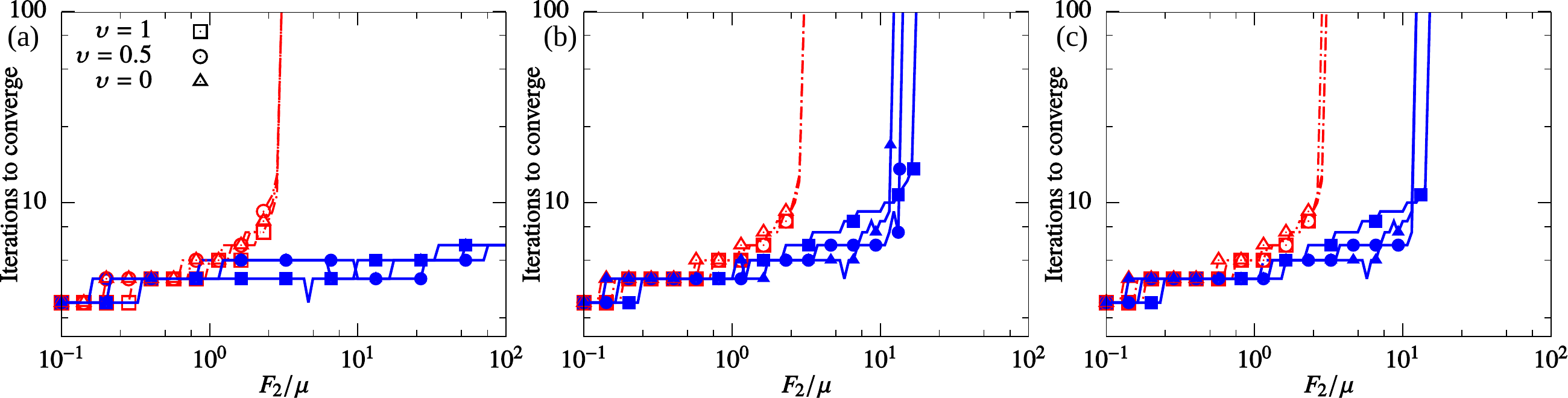}
    \caption{Convergence plots of the cases of geometric nonlinearity for the Mooney-Rivlin model using the standard formulation (red dashed line and open symbols) and arctan formulation (blue solid line and filled symbols) with varying $\upsilon$: iterations taken to converge for the (a) 1D volume-preserving compression, (b) axisymmetric compression ($K/\mu=10$), and (c) 3D uniaxial compression  ($K/\mu=10$) problems. \label{1D-spring-arctan-fig}}
\end{figure}

To test the arctan formulation, we again start with the uniaxial axisymmetric compression under volume-preserving constraint. Using the standard method, the number of iterations required to converge increases with the applied load (Fig.~\ref{1D-spring-arctan-fig}a), and the method fails to converge for $F_2/\mu>3$. For small load steps, the arctan formulation performs similar to the standard formulation and converges in the same number of iterations. However, using the arctan formulation the number of iterations does not increase with an increasing load step and the formulation converges in less than ten iterations for load steps as large as $F_2/\mu=100$. For both formulations, the iterations required to converge are largely insensitive to the parameter $\upsilon$. Varying $F_1/\mu$ for the neo-Hookean material does not make any difference in its convergence behavior (result figure skipped for brevity). 

\subsubsection{Axisymmetric compression}

When we relax the incompressibility constraint ($K/\mu=10$), the overall trend remains the same and the arctan method outperforms the standard method (Fig.~\ref{1D-spring-arctan-fig}b). The standard formulation performs exactly the same as the previous case. However, the improvement in convergence using the arctan formulation decreases slightly, and the number of iterations required to converge slowly increases with an increasing load step before failing to converge at $F_2/\mu\approx15$. 

\subsubsection{3D uniaxial compression}
\label{3D uniaxial compression}
\begin{figure}[!tbp]
  \centering
  \includegraphics[width=0.9\textwidth]{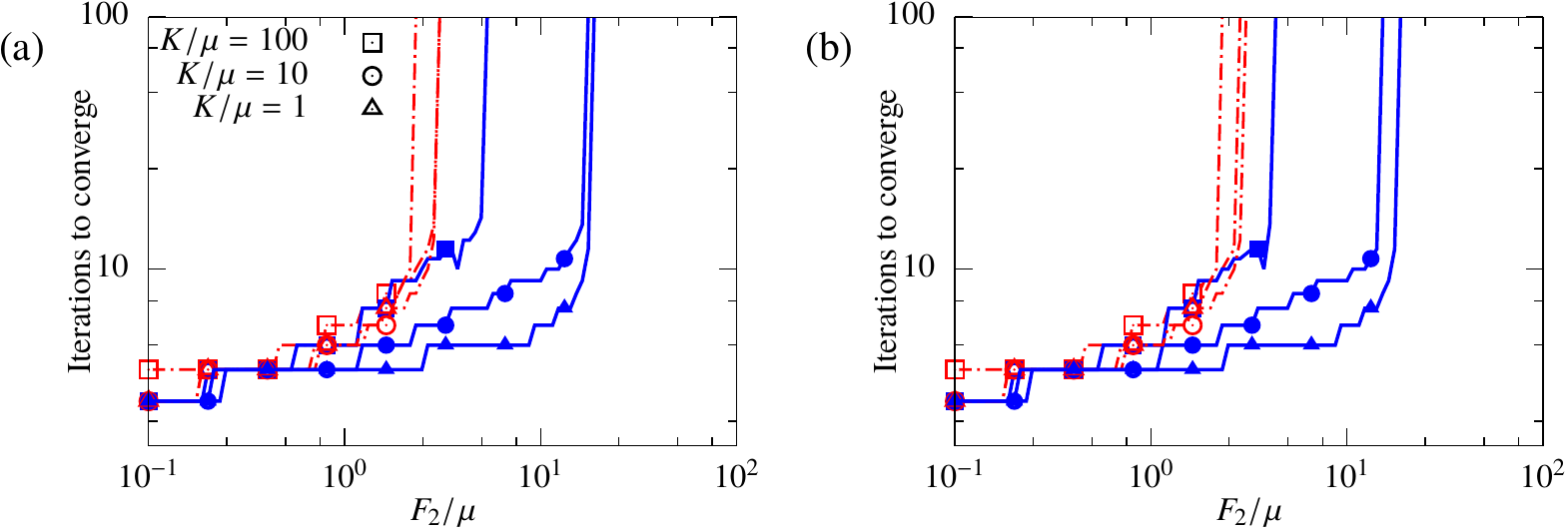}
    \caption{Effect of compressibility on the convergence in the cases of the neo-Hookean model using the standard formulation (red dashed line and open symbols) and arctan formulation (blue solid line and filled symbols) with (a) axisymmetric assumption and (b) general 3D case. \label{arctan-bulk}}
\end{figure}

In the case of general three-dimensional uniaxial compression, both formulations perform similar to the axisymmetric case (Fig.\ref{1D-spring-arctan-fig}c). In spite of more degrees of freedom, there is no noticeable change in the convergence behavior. Overall, it is clear that, for compression cases, the arctan formulation allows significantly larger load steps compared to the standard formulation. However, compared to the log formulation, changing the bulk modulus has an opposite effect on the arctan formulation. The improvement in convergence behavior due to the arctan framework is reduced as the bulk modulus is increased (Fig.~\ref{arctan-bulk}). This effect is observed for both axisymmetric and 3D problems.

\section{Practical examples}
\label{practical-problems}

After studying the convergence of the proposed formulation using relatively simple numerical examples, we turn to solve three practical problems and showcase the advantage of using the proposed method. These problems are motivated by realistic nonlinear situations: pressurization of the aorta, indentation experiment of the elastic solid, and biaxial stretch experiment of the planar solid. For each of them, we compare the convergence behavior by plotting the $L_2$ norm of the displacement error. The displacement error is defined as $\left\| {{\bf{U}} - {{\bf{U}}^i}} \right\|/\left\| {\bf{U}} \right\| $, where ${\bf{U}}^i$ and $\bf{U}$ are the nodal displacement vectors at the $i$-th iteration and final iteration, respectively, on a chosen subset of nodes.

\begin{figure}[!h]
  \centering
    \includegraphics[width=1.0\textwidth]{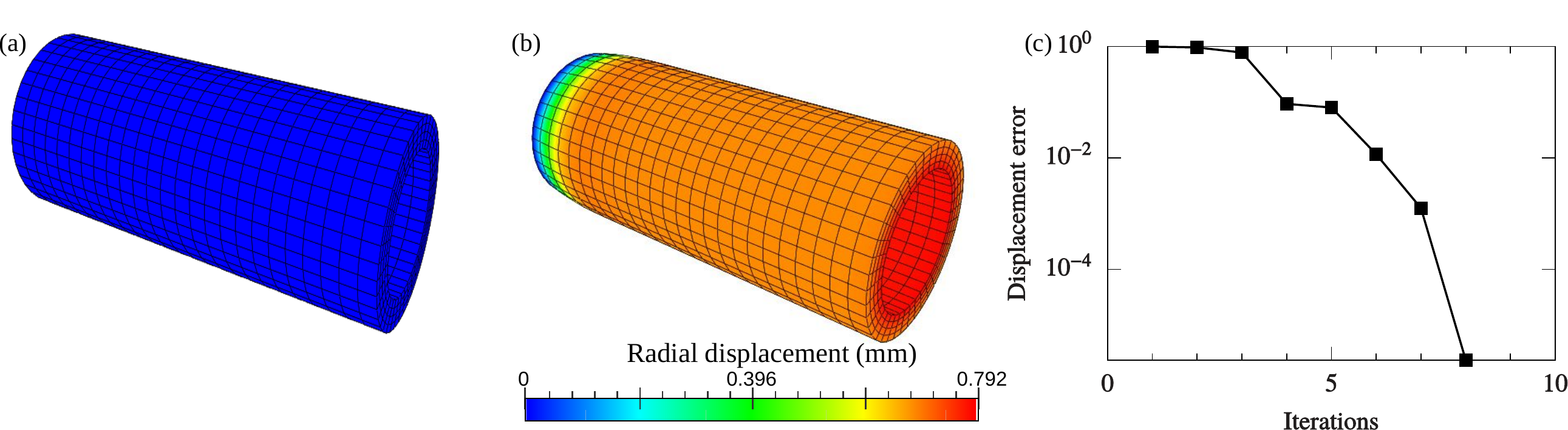} 
    \caption{Pressurization simulation: (a) mesh of the a thin wall tube; (b) simulation result using the log formulation: the deformed configuration colored by the resulting radial displacement at an internal pressure of 0.2 kPa; (c) displacement error norm versus iteration using log formulation for the second load step. \label{tube-convergence}}
\end{figure} 

\subsection{Pressurization}
Biomechanical characterization of aorta is an important topic in cardiovascular research. We model the aorta as a 5 cm long thick-walled tube with uniform inner and outer radii of $0.7$ cm and $1$ cm, respectively (Fig.~\ref{tube-convergence}). The aortic tissue is modeled using an isotropic VW constitutive law, with the nonlinear parameter $B=50$, stiffness $A=0.5$ kPa, and bulk modulus $K=10$ kPa. All degrees of freedom at the left end of the tube are fixed, while the right end is free. The tube is discretized into 4000 uniform hexahedral finite elements with 4, 40, and 25 elements along radial, circumferential, and axial directions of the tube, respectively (Fig.~\ref{tube-convergence}a). The element sizes are chosen to obtain a fine enough discretization while keeping the computational expense reasonable. Aorta undergoes a periodic pressurization during the cardiac cycle, which is modeled as a normal static uniform pressure load on the internal surface. Unlike the uniaxial extension problem used in the previous section, here the applied load is not aligned with a particular axis of deformation. The solution is first computed at $2\times 10^{-5}$ kPa of pressure (first load step), and then the pressure is increased to $0.2$ kPa (second load step). This pressure induces roughly $8\%$ increase of the outer radius (Fig.~\ref{tube-convergence}b). As the nonlinear parameter $B$ is large, this problem is highly nonlinear and numerically challenging to solve. We observe that the log formulation yields a converged solution in eight iterations (Fig.~\ref{tube-convergence}c), whereas the standard formulation fails to converge.

\begin{figure}[!h]
  \centering
  \includegraphics[width=\textwidth]{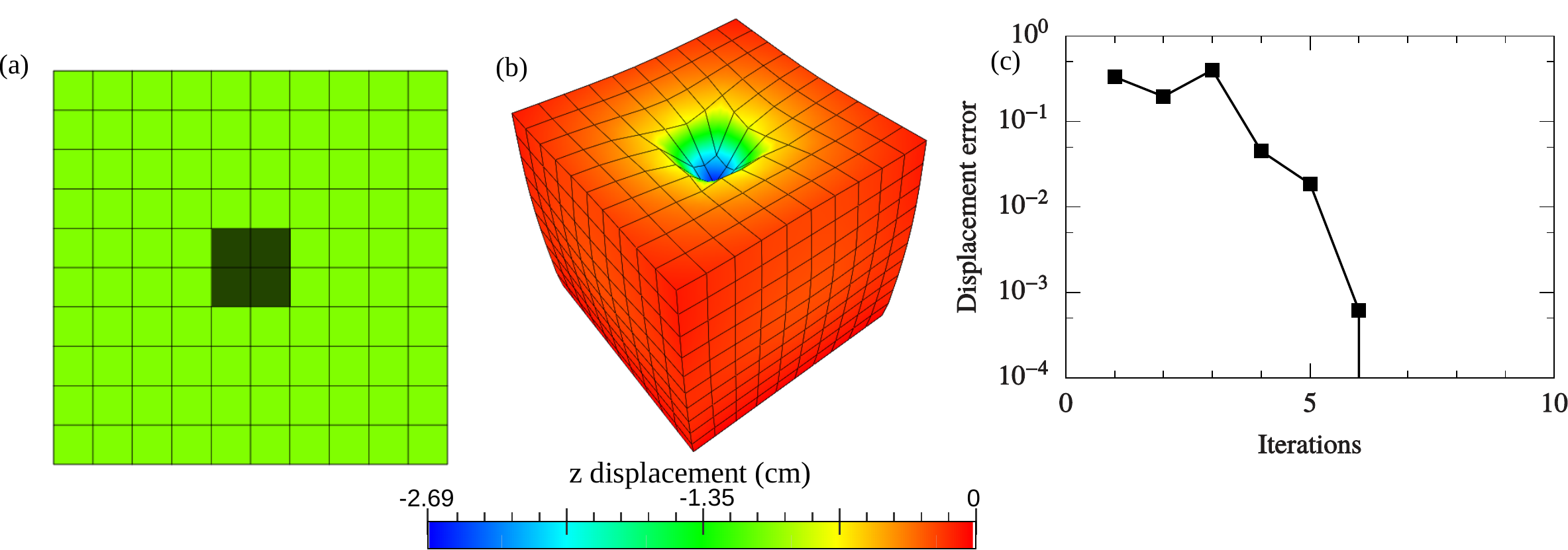}
  \caption{Indentation simulation: (a) the top surface of the cubic sample where we apply the local uniform force (on the shaded region); (b) Resulting deformed configuration using the arctan formulation when the total applied force is 360 N; (c) norm of the displacement error versus iteration using arctan formulation. \label{3D-VM-local-fig}}
\end{figure} 

\subsection{Indentation}
Indentation is frequently used to determine the local mechanical properties of solids. To simulate an indentation experiment, we apply local compressive pressure on a small central region at the top surface of a $10\times10\times10$ cm$^3$ cubic solid (Fig.~\ref{3D-VM-local-fig}a) and fix all degrees of freedom at the bottom surface. Unlike the previous uniaxial compression, this setup induces a non-uniform deformation (Fig.~\ref{3D-VM-local-fig}b). We use the neo-Hookean model with shear and bulk moduli $\mu=0.2$ MPa and $K=1$ MPa, respectively. Similar to Section \ref{results}, we use two load steps: $F_1=0.036$ N and $F_2=360$ N (corresponding to a local pressure of 0.9 MPa). This local pressure leads to a large local deformation with stretches up to $27\%$ (Fig.~\ref{3D-VM-local-fig}b).  We observe that the standard method fails to converge, whereas the arctan formulation takes only seven iterations to converge (Fig.~\ref{3D-VM-local-fig}c).  

\begin{figure}[!h]
\centering
\includegraphics[width=\textwidth]{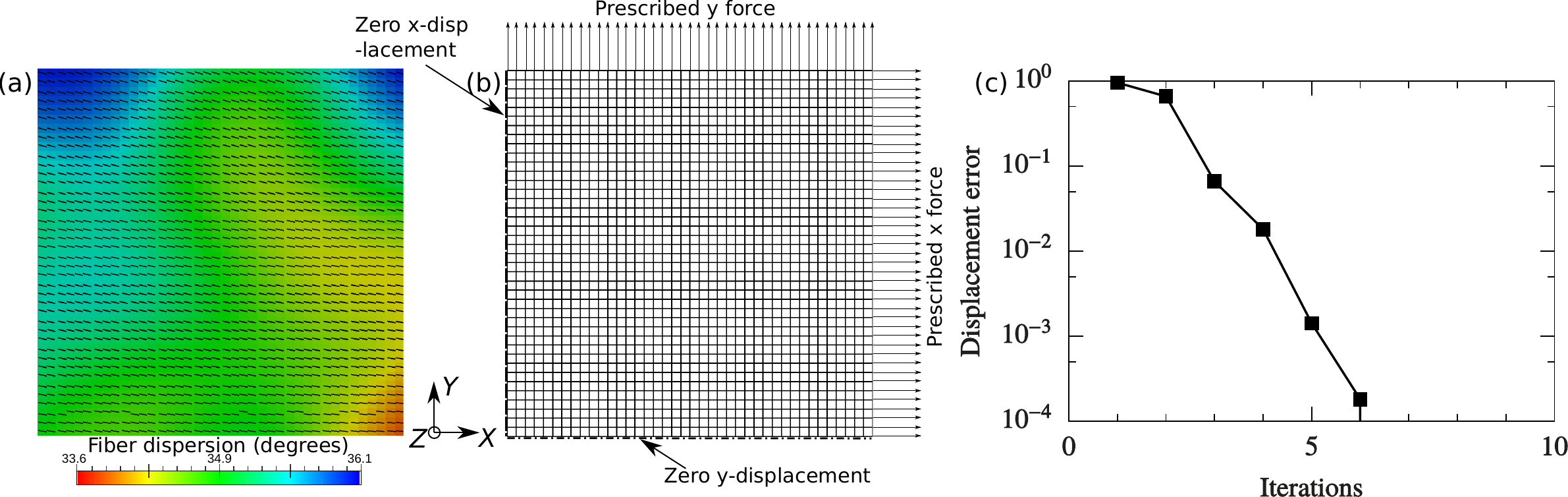}
\caption{Biaxial testing simulation: (a) microstructure of an anisotropic thin planar tissue with fibers aligned predominantly in the x-direction; (b) simulation setup of its biaxial testing; (c) norm of displacement error versus iteration using log formulation. \label{2D-SM-fig}}
\end{figure}  

\subsection{Biaxial testing}
Lastly, we consider the simulation of the biaxial testing of a thin tissue sample for a different constitutive model with an exponential function: the simplified structural model \citep{fan2014simulation}. This model also produces a stress-strain relation $\sigma \sim A \exp(B\epsilon)$, with $A=0.02$ MPa and $B=44.6$ resulting in a highly nonlinear response. Importantly, the simplified structural model results in an anisotropic response. Other parameters related to the fiber orientation function and ground matrix are the same as those used in our previous study \citep{aggarwal2017improved}. A $3.3\times3.3$ mm$^2$ thin planar tissue sample is subjected to biaxial stretch (Fig.~\ref{2D-SM-fig}). We restrict the axial motion of the left and bottom edges, apply uniform nodal forces on the right and top edges, and use two load steps to solve this problem. The nodal forces applied on the right and top edges at the the first load step are $4.08\times10^{-5}$ N and $2.57\times10^{-5}$ N, respectively. Both nodal forces are increased by 5000 times at the second load step. This induces $12\%$ and $15\%$ extension in x- and y-directions, respectively. When the exponent parameter $B$ is large, this problem becomes highly nonlinear. We observe that the standard method is nonconvergent, while it only takes seven iterations for the new method to converge (Fig.~\ref{2D-SM-fig}c).

\section{Discussion}
\label{Discussion}
\subsection{Linear vs. nonlinear problems}
The initial development of the field of computational engineering was based on linear systems, and most of the techniques were developed for linear equations. These developments included numerical techniques in linear algebra, such as preconditioning, which are widely used in all computations \cite{greenbaum1997iterative}. As the scope of computational methods expanded, their application to problems with high nonlinearity became more commonplace. 

The most common approach is to linearize the governing nonlinear equations, so that the same linear algebra tools can be employed. This is in spite of the fact that nonlinear problems demonstrate significant differences compared to linear problems, especially those related to convergence and stability of the numerical methods. Highly nonlinear problems are notoriously hard to converge and can even lead to unstable numerical schemes that are unconditionally stable for linear problems. A commonly used solution is to divide the problem into several small steps, such that the current guess is always close to the solution and, thus, linearization is applicable.

\subsection{Significance of presented framework}

In this paper, we present a novel viewpoint by transforming the governing equations \emph{before} linearization. The transformation is determined so as to decrease the nonlinearity of the resulting space-discretized problem \eqref{eq-tau}. After the nonlinearity has been reduced in the transformed formulation, the linearization is expected to be valid in a larger neighborhood of the solution. Hence, the problem will be able to converge even for larger steps, leading to a significantly faster solution. It is important to note that the rate of convergence in the neighborhood of solution remains the same using the proposed formulation. Instead, the improvement stems from our ability to take larger load steps.

In this study, we focused on static elasticity equations, so that the governing equations can be written in a specific form \eqref{original-eq}. After applying the transformation and linearization, we found that the new equation \eqref{transformed-linear-eq} has a very similar structure to the standard equations \eqref{linear-eq}. In fact, an assumption on the geometric stiffness matrix gave us a formulation, where the stiffness matrix remains the same and only the residual vector is changed \eqref{transformed-linear-eq3}. This relatively small change makes our formulation incredibly easy to be implemented in an existing solver. 

It is worthwhile discussing a similarity between the presented framework and the idea of preconditioning. Most of the previous work has been done on linear preconditioners, where a linearized system of equations $\mathbf{A}\boldsymbol{x}=\boldsymbol{b}$ is transformed by pre-multiplying by a matrix so as to improve the numerical conditioning of the problem. Instead of working on the linearized system, the presented framework targets the nonlinear equations before linearization. However, the underlying idea of transforming the equations for improving the numerical solution has a parallel, making the proposed framework similar to \emph{nonlinear} preconditioning --- a relatively unexplored field \cite{cai2002nonlinearly}. Furthermore, usually preconditioners are designed based on purely mathematical properties of the problem. Instead, here we propose to design a transformation based on the physical characteristics of the problem and the predominant type of nonlinearity it contains.

\subsection{Presented cases}

One key aspect of the proposed framework is that a transformation must be determined for a \emph{known} nonlinearity. As the focus is on highly nonlinear problems, we picked two distinctly different nonlinearities --- exponential-type material nonlinearity and geometric nonlinearity in large compression. We determined that a log transformation decreases the exponential-type nonlinearity, whereas an arctan formulation decreases the geometric nonlinearity in compression. Both formulations require the first load step to be solved using the standard formulation. This is because a non-zero internal force is required to apply log transformation and the scaling factor $\alpha_i$ needs to be determined before applying the arctan formulation. Therefore, the transformed equations were used only in the second load step. 

We first solved a simple uniaxial extension/compression problem (Section~\ref{numerical-examples}) in 3D and its reduction under axisymmetric and incompressibility constraints. Overall, the proposed formulations led to an improved convergence, although to different extents. The improvement delivered by using the log formulation was dependent on the problem nonlinearity, which depends on the exponent parameter $B$. Interestingly, the performance of the log formulation improved further for solids with a larger bulk modulus. For the arctan formulation, the improvement did not change with material parameters. This is expected since the nonlinearity in this case does not originate from the material properties. However, the effect of bulk modulus was opposite, and the formulation worked better for solids with a smaller bulk modulus.

We also tested the proposed formulations for three practical problems -- indentation testing, biaxial testing and pressurization of an arterial tube (Section~\ref{practical-problems}). These problems presented additional complexities, such as a non-uniform deformation and external force not being aligned with a coordinate axis. Furthermore, an anisotropic material model was used in the biaxial testing simulation, although with the exponential nonlinearity as well. We were able to solve all three problems in significantly less computational time. 

\begin{figure}[!tbp]
  \centering
\includegraphics[width=0.9\textwidth]{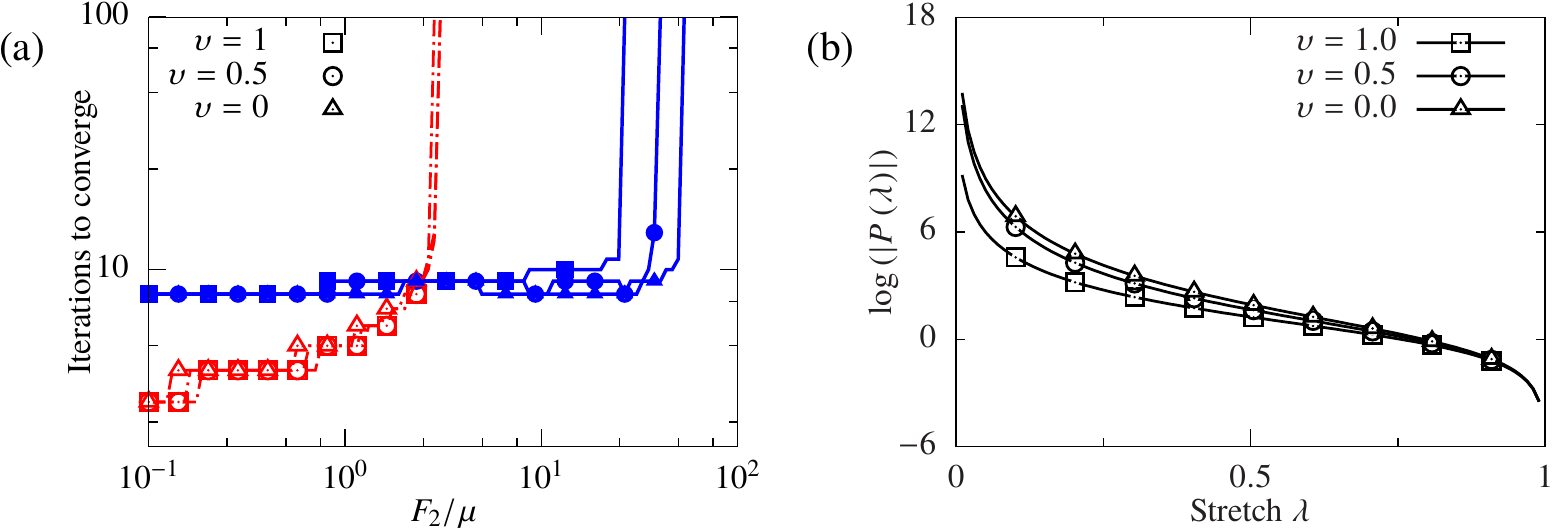}
    \caption{(a) Iterations taken to converge for the 3D uniaxial compression problem with Mooney-Rivlin model using the standard formulation (red dashed line and open symbols) and log formulation (blue solid line and filled symbols) for varying $\upsilon$; (b) the transformed stress under compression using log transformation \label{3D-log-MM-fig}}
\end{figure}

We recognize that transformation that decreases a given nonlinearity may not be unique. For instance, Fig~\ref{3D-log-MM-fig}a presents the iterations taken for the 3D uniaxial compression problem using log formation instead of the arctan formation. We find that the log formation also improves the convergence at large load steps, but takes more iterations to converge at small load steps. In order to understand this surprising result, we compare the arctan and logarithm of 1st PK stress under compression (Figs.~\ref{mooney-fig}b, \ref{3D-log-MM-fig}b). We note that the logarithm makes the stress-stretch relation roughly linear for stretch ($\lambda$) between 0.1 and 0.9, while it increases the nonlinearity at small loads (i.e. stretch close to 1). This is why the log formulation improves the convergence for the compression case for large loads, but it requires more iterations to converge at small loads. The improved convergence of the log formulation even for the compression case is an unexpected result, thus making this formulation an even more attractive choice for soft tissues. 

\subsection{Limitations and future work}

In this study, transformations were determined for two specific cases, and this work needs to be extended to other types of nonlinearities. Here we focused on problems where the nonlinearity originates from a single recognizable source. The next important step will be to determine how this approach can be applied to cases where different nonlinearities are coupled together. We also presented a few different choices in the implementation (\ref{transformed-linear-eq}-\ref{transformed-linear-eq3}), but only the one with the simplest modification was numerically tested.  Although we presented results using Newton's method, the reduction in nonlinearity means that a similar improvement is expected using other popular solvers such as quasi-Newton methods \cite{zhu1997algorithm}.

All of the problems tested here were driven by traction force boundary condition. If the deformation is driven by applied displacement at the boundary, the application of the presented formulation will depend on how the fixed DOFs on Dirichlet boundary nodes are treated. If the fixed DOFs are eliminated before linearization, the equations reduce to our standard form \eqref{original-eq} with the external force derived from the applied displacements. In such an implementation, our proposed formulation will be applicable as is. However, for other implementations of boundary conditions, such as using Lagrange multiplier or elimination after linearization, the formulation will need modifications.

In addition to working on other types of nonlinearities and addressing the limitations described above, in the future we will determine a unified approach to decrease the nonlinearity of a general problem. Since it may not be possible to determine a simple transformation $\mathcal{T}$ for every case, a numerical approach may be more suitable. Here we only presented a framework for static elasticity problems. In the future, we will extend this framework to implicit and explicit dynamic cases as well. Besides, the proposed framework will be applied to other nonlinear problems in mechanics, such as
inelasticity. Importantly, this proposed framework can be used to motivate design of novel nonlinear preconditioners.

\section{Conclusion}
\label{Conclusion}
In this paper, we presented a framework to improve the convergence behavior of FE-based numerical schemes used for solving nonlinear elasticity problems in the static case. This approach requires recognizing the main source of nonlinearity for a specific problem and determining a transformation that reduces its nonlinearity. We tested the feasibility of our method in two scenarios: material non-linearity induced by the exponential function and geometric non-linearity due to the compression. We determined two different transformations for these problems and observed that the proposed method significantly increases the permissible load step. This novel framework is simple to implement and can be easily integrated into any existing finite element solver. Hence, this approach has the potential of addressing the convergence issues existing in many computational techniques for highly nonlinear elasticity problems.



\section*{Acknowledgement}
This work was supported by Welsh Government and Higher Education Funding Council for Wales through the S\^{e}r Cymru National Research Network in Advanced Engineering and Materials (Grant No. F28), and the Engineering and Physical Sciences Research Council of the UK (Grant No. EP/P018912/1).

\section*{References}
\bibliographystyle{elsarticle-num} 
\bibliography{refs}


\end{document}